\documentclass[a4paper, 11pt]{article}

\usepackage[T2A]{fontenc} 
\usepackage[cp1251]{inputenc} 

\usepackage[english]{babel}
\usepackage[tbtags]{amsmath}
\usepackage{amsfonts,amssymb,mathrsfs,amscd}

\usepackage{misccorr} 
\usepackage{indentfirst} 

\usepackage{graphicx}

\newtheorem{theorem}{Theorem}[section]
\newtheorem{lemma}{Lemma}[section]
\newtheorem{definition}{Definition}[section]
\newtheorem{remark}{Remark}[section]

\newtheorem{corollary}{Corollary}[section]
\newtheorem{proposition}{Proposition}[section]

\usepackage{xcolor}
\usepackage{hyperref}
 
\definecolor{linkcolor}{HTML}{799B03} 
\definecolor{urlcolor}{HTML}{799B03} 
 
\hypersetup{pdfstartview=FitH,  linkcolor=linkcolor,urlcolor=urlcolor, colorlinks=true}
\begin{document}


\title{General Markov Chains:  Cycles of Finitely Additive Measures
          and Сlassical Cycles of States }

\author{Alexander ~I.~Zhdanok \footnote{Institute for Information Transmission Problems of the Russian Academy of Science, Moscow, Russia;} \footnote{Tuvinian Institute for Exploration of Natural Resources of the Siberian Branch RAS, Kyzyl, Republic of Tuva, Russia;} 
\footnote{This work was supported by the Russian Foundation for Basic Research (project No. 20-01-00575-a)}\\
{\it e-mail: zhdanok@inbox.ru}
}

\date{} 

\maketitle

{\bf Abstract.}
General Markov chains in an arbitrary phase space are considered in the framework of the operator treatment. Markov operators continue from the space of countably additive measures to the space of finitely additive measures. Cycles of measures generated by the corresponding operator are constructed, and algebraic operations on them are introduced. One of the main results obtained is that any cycle of finitely additive measures can be uniquely decomposed into the coordinate-wise sum of a cycle of countably additive measures and a cycle of purely finitely additive measures.We have proved theorems on the conditions and consequences of consistency cycles of measures with cycles of sets of states of General Markov chains. A theorem is proved (under certain conditions) that if a finitely additive cycle of a Markov chain is unique, then it is countably additive.  
\\




\pagestyle{plain} 

\newpage
\section{Introduction}
\label{S:1}

The considered general Markov chains (MC) are random processes with an arbitrary phase space, with discrete time, and homogeneous in time. MCs are given by the usual transition probability, countably additive in the second argument, which generates two Markov operators $T$ and $A$ in the space of measurable functions and in the space of countably additive measures, respectively. Thus, we use the operator treatment in the theory of general MCs, proposed in 1937 by N. Kryloff and N. Bogolyuboff, and developed in detail in the article \cite{YoKa:1}. Later, in a number of works by different authors, an extension of the Markov operator $A$ to the space of finitely additive measures was carried out, which turned $A$ into an operator topologically conjugate to the operator $T$, and opened up new possibilities in the development of the operator treatment. Within the framework of such a scheme, we carry out here the study of cycles of measures of general MC. In this case, we use a number of information on the general theory of finitely additive measures from the sources ~\cite{YoHew1} and ~\cite{DS1}.

In the ergodic theory of MC, one usually distinguishes in the space of its states ergodic classes and their cyclic subclasses, if such exist (see, for example, ~\cite{Rev1}). However, in the general phase space, the study of such sets has its natural limitations. Therefore, in some cases it is more convenient to use not cycles of sets, but cycles of measures generated by the Markov operator A.

In this paper we propose a corresponding construction for cyclic finitely additive measures of MC on an arbitrary measurable space.
We also study cycles of countably additive and purely finitely additive measures, and their relationship with classical cycles of sets of States of an arbitrary Markov chain.
 In the proof of the theorems presented here, we also use some results of papers ~\cite{Zhd01} and ~\cite{Zhd02}.

\section{Finitely additive Measures and Markov operators}
\label{S:2}

Let $X$ be an arbitrary infinite set and $\Sigma$ the sigma-algebra of its subsets containing all one-point subsets from $X$. Let $B(X,\Sigma)$ denote the Banach space of bounded $\Sigma$ -measurable functions $f: X\to R$ with sup-norm. 

We also consider Banach spaces of bounded measures $\mu : \Sigma\to R$, with the norm equal to the total variation of the measure $\mu$ (but you can also use the equivalent sup-norm): 

$ba(X,\Sigma)$ is the space of finitely additive measures, 

$ca(X,\Sigma)$ is the space of countably additive measures. 

If $\mu\ge{0}$, then $||\mu||=\mu(X)$.


\begin{definition} (see \cite{YoHew1}). \label{D:2.1}
A finitely additive nonnegative measure $\mu$ is called purely finitely additive (pure charge, pure mean) if any countably additive measure $\lambda$ satisfying the condition $0\le\lambda\le\mu$ is identically zero. An alternating measure $\mu$ is called purely finitely additive if both components of its Jordan decomposition are purely finitely additive.
\end{definition}

Any finitely additive measure $\mu$ can be uniquely expanded into the sum  $\mu =\mu_1 + \mu_2$, where $\mu_1$ is countably additive and $\mu_2$ is a purely finitely additive measure (the Alexandroff  -Yosida-Hewitt decomposition, see \cite{YoHew1} and \cite{DS1}).

Purely finitely additive measures also form a Banach space $pfa(X,\Sigma)$ with the same norm, $ba(X,\Sigma) =ca(X,\Sigma) \oplus pfa(X,\Sigma)$.\\

{\bf Examples 2.1}. Here are two examples of purely finitely additive measures.

Let $X=[0, 1] \subset \mathbb R$ ($\mathbb R = (- \infty; + \infty)$) and $\Sigma = \mathfrak B$ ( $\mathfrak B$ - Borel sigma algebra).
There is (proved) a finitely additive measure $\mu \colon \mathfrak B \to \mathbb R$, $\mu \in S_{ba}$, such that for any $ \varepsilon> 0 $ the following holds:
  $$\mu((0, \varepsilon)) = 1, \;\; \mu([\varepsilon, 1]) = 0, \;\; \mu(\{0\}) = 0.$$
 We can say that the measure $ \mu $ fixes the unit mass arbitrarily close to zero (on the right), but not at zero. According to \cite{YoHew1}, such a measure is purely finitely additive, but it is not the only one. It is known that the cardinality of a family of such measures located ``near zero (on the right)'' is not less than $ 2^{2^{\aleph_{0}}} = 2^{c} $ (hypercontinuum). And the same family of purely finitely additive measures exists ``near each point $ x_{0} \in [0, 1] $ (to the right, or to the left, or both there, and there)''.\\  
 
 {\bf Examples 2.2}.

Let $X=\mathbb R = (- \infty; + \infty)$ and $\Sigma = \mathfrak B$.
There is (proved) a finitely additive measure $\mu \colon \mathfrak B \to \mathbb R$, $\mu \in S_{ba}$, such that for any $x \in \mathbb R$ the following holds:
$$\mu((x, \infty)) = 1, \;\; \mu((-\infty, x)) = 0, \;\; \mu(\{x\}) = 0.$$
We can say that the measure $ \mu $ fixes the unit mass arbitrarily far, ``near infinity''. This measure is also purely finitely additive. And there are also a lot of such measures.
\\

We denote the sets of measures:

$S_{ba}=\{\mu\in{ba(X,\Sigma)}:\mu\ge{0}, ||\mu||=1\},$ $S_{ca}=\{\mu\in{ca(X,\Sigma)}:\mu\ge{0}, ||\mu||=1\},$

$S_{pfa}=\{\mu\in{pfa(X,\Sigma)}:\mu\ge{0}, ||\mu||=1\}.$

All measures from these sets will be called probabilistic. 

Markov chains (MC) on a measurable space $(X,\Sigma)$ are given by their transition function (probability) $p(x,E), x\in X, E\in\Sigma$, under the usual conditions: 

\begin{enumerate}
\item $0\le p(x,E) \le{1}, p(x,X)=1, \forall{x}\in{X}, \forall{E}\in\Sigma$;
\item $p(\cdot,E)\in{B(X,\Sigma)}, \forall{E}\in\Sigma$;
\item $p(x,\cdot)\in{ca(X,\Sigma)}, \forall{x}\in{X}$.
\end{enumerate}

We emphasize that our transition function is a countably additive measure in the second argument, i.e. we consider classical MCs. 

The transition function generates two Markov linear bounded positive integral operators: 
 
$T: B(X,\Sigma)\to{B(X,\Sigma)}, (Tf)(x)=Tf(x)=\int\limits_X f(y)p(x,dy),$

$\forall{f\in{B(X,\Sigma)}},\forall{x}\in{X};$ \medskip

$A: ca(X,\Sigma)\to{ca(X,\Sigma)}, (A\mu)(E)= A\mu(E)=\int\limits_X p(x,E)\mu(dx),$

$\forall{\mu\in{ca(X,\Sigma)}},\forall{E}\in\Sigma.$ \medskip

Let the initial measure be $\mu_0\in{S_{ca}}$. Then the iterative sequence of countably additive probability measures $\mu_{n}=A\mu_{n-1}\in{S_{ca}},n\in{N}$, is usually identified with the Markov chain. 

Topologically conjugate to the space $B(X,\Sigma)$ is the (isomorphic) space of finitely additive measures: $B^*(X,\Sigma)=ba(X,\Sigma)$ (see, for example, \cite{DS1}). Moreover, the operator $T^*:ba(X,\Sigma)\to{ba(X,\Sigma)}$ is topologically conjugate to the operator $T$:
$$
T^*\mu(E)=\int\limits_X p(x,E)\mu(dx), \forall{\mu\in{ba(X,\Sigma)}}, \forall E\in\Sigma.
$$ 

The operator $T^*$ is the only bounded continuation of the operator $A$ to the entire space $ba(X,\Sigma)$ while preserving its analytic form. The operator $T^*$ has its own invariant subspace $ca(X,\Sigma)$, i.e. $T^*[ca(X,\Sigma)] \subset{ca(X,\Sigma)}$, on which it matches the original operator $A$. The construction of the Markov operators $T$ and $T^*$ is now functionally closed. We shall continue to denote the operator $T^*$ as $A$.

In such a setting, it is natural to admit to consideration also the Markov sequences of probabilistic finitely additive measures $ \mu_0\in{S_{ba}}, \mu_n = A\mu_{n-1}\in{S_{ba}}, n\in{N}$, keeping the countable additivity of the transition function $p(x,\cdot)$ with respect to the second argument.

\section{Cycles of Measures and their properties} 
\label{S:3}


\begin{definition}
If $A\mu = \mu$ holds for some positive finitely additive measure $\mu$, then we call such a measure {\it invariant} for the operator $A$ (and for the Markov chain).
\end{definition}

We denote the sets of all probability invariant measures for the  operator $A$:

$\Delta_{ba}=\{\mu\in{S_{ba}}: \mu=A\mu\}$,

 $\Delta_{ca}=\{\mu\in{S_{ca}}: \mu=A\mu\}$, $\Delta_{pfa}=\{\mu\in{S_{pfa}}: \mu=A\mu\}$.

A classical countably additive Markov chain may or may not have invariant countably additive probability measures, i. e. possibly $\Delta_{ca}=\O$ (for example, for a symmetric walk on $Z$).

In \cite[Theorem 2.2]{Si1} \v{S}idak proved that any countably additive MC on an arbitrary measurable space $(X,\Sigma)$ extended to the space of finitely additive measures has at least one invariant finitely additive measure, i.e. always $\Delta_{ba}\ne\O$.
\v{S}idak in \cite[Theorem 2.5]{Si1} also established in the general case that if a finitely additive measure $\mu$ is invariant $A\mu=\mu$, and $\mu=\mu_1+\mu_2$ is its decomposition into are countably additive and purely finitely additive components, then each of them is also invariant: $A\mu_1=\mu_1$, $A\mu_2=\mu_2$. Therefore,  it suffices to study invariant measures from $\Delta_{ca}$ and from $\Delta_{pfa}$, separately.

The study of questions of the existence and properties of invariant finitely additive measures of general Markov chains is also carried out in the author's preprints ~\cite{Zhd10}, ~\cite{Zhd12}, ~\cite{Zhd14}.

 \begin{definition} \label{D:3.2}
 A finite numbered set of pairwise different positive finitely additive measures $K=\{ \mu_{1}$, $\mu_{2}$, ... , $\mu_{m} \}$ will be called a {\it cycle measures of an operator $A$} of a given Markov chain (or a cycle of measures MC ) if
$$
A\mu_{1}=\mu_{2}, A\mu_{2}=\mu_{3}, ... , A\mu_{m-1}=\mu_{m}, A\mu_{m}=\mu_{1}.
$$
\end{definition}

Such cycles will be called finitely additive. 
The number $m\ge 1$ will be called the {\it cycle} {\it period}, and the measures $\mu_{1}, \mu_{2}, ... , \mu_{m}$ -- {\it cyclic} {\it measures}. Unnormalized cycles will also be used below.

If $K=\{ \mu_{1}, \mu_{2}, ... , \mu_{m} \}$ is a MC cycle, then, obviously,
$$
A^{m}{\mu_{1}}=\mu_{1}, A^{m}{\mu_{2}}=\mu_{2}, ... , A^{m}{\mu_{m}}=\mu_{m},
$$
i.e. all cyclic measures $\mu_{i}$ are invariant for the operator $A^{m}$ and $ A^{m}(K)=K$.

The following well-known statement is obviousю.
  Let $K=\{ \mu_{1}, \mu_{2}, ..., \mu_{m} \} $ be a cycle of finitely additive measures. Then the measure
$$
\mu = \frac{1}{m} \sum_{k=1}^{m}\mu_{k} = \frac{1}{m} \sum_{k=1}^{m}A^{k-1}\mu_{1}
$$
is invariant for the operator $A$, i.e. $A\mu = \mu$ (here $ A ^ {0} $ is the identity operator).


\begin{definition} 
The measure constructed above will be called the {\it mean cycle measure} $K$.
\end{definition}

 
 \begin{definition} 
We call each method of choosing a measure $\mu_ {1} $ in $ K $ an operation renumbering a cycle $ K $.
\end{definition}

 \begin{definition} 
   We say that two cycles of the same period $K^1$ and $K^2$are identical if there is a renumbering of cycles $ K ^ 1 $ or $ K ^ 2 $ such that all their cyclic measures with the same numbers match. In this case, we will write $ K ^ 1 = K ^ 2 $. Instead of the words `` identical cycles '', we will still say the words `` equal cycles ''.
\end{definition}

Obviously, for the cycles to be equal, it is sufficient that their first measures coincide.

Hereinafter, it is convenient to call cyclic measures $ \mu_ {i}, i = 1, ..., m, $ {\it cycle coordinates} $ K $. 


\begin{definition}
{\it By the operation of multiplying a cycle of measures $K=\{ \mu_{1}, \mu_{2}, ..., \mu_{m} \}$ by a number $\gamma > 0$} we mean the construction of a cycle of measures $\gamma K=\{ \gamma \mu_{1}, \gamma \mu_{2}, ..., \gamma \mu_{m} \}$.
\end{definition}

Since the operator $ A $ is isometric in the cone of positive measures, all cyclic measures of one cycle $K=\{ \mu_{1}$, $\mu_{2}$, ... , $\mu_{m} \}$ have the same norm $\| \mu_{1}\|=\|\mu_{2}\|= , ... , = \| \mu_{m}\| = \|\mu\| $, which is naturally called the norm $ \| K \| $ of the cycle $ K $ itself.

To give the cycle a probabilistic meaning, it is sufficient to multiply it coordinatewise by the normalizing factor $\gamma = \frac{1}{\|\mu\|}$: $\hat K = \gamma \cdot K = \{ \gamma \mu_{1}$, $\gamma \mu_{2}$, ... , $\gamma \mu_{m} \}$. We obtain a probability cycle with the norm $\|\hat {K}\|=1, \hat {K} \subset S_{ba}$.


\begin{definition}
Let there be given two cycles of measures of the same MC $K^{1}=\{ \mu^{1}_{1}, \mu^{1}_{2}, ..., \mu^{1}_{m} \}$  and   
$K^{2}=\{ \mu^{2}_{1}, \mu^{2}_{2}, ..., \mu^{2}_{m} \}$ of the same period $m$. We call the {\it sum of cycles} $K^{1}$ and $K^{2}$ the following set of measures $K=K^{1}+K^{2}=\{ \mu^{1}_{1}+\mu^{2}_{1}, ..., \mu^{1}_{m}+\mu^{2}_{m} \}$ derived from $K^{1}$ and $K^{2}$ coordinatewise addition.
\end{definition}

The measure spaces are semi-ordered by the natural order relation.
In them one can introduce the notion of infimum $inf \{\mu_{1} , \mu_{2}\}=\mu_{1} \wedge \mu_{2}$ and supremum
$sup \{\mu_{1} , \mu_{2}\}=\mu_{1} \vee \mu_{2}$, which are also contained in these spaces.
Thus, the measure spaces $ba(X, \Sigma), ca(X, \Sigma)$ and $pfa(X, \Sigma)$ are lattices ($K$-lineals).

In~\cite {YoHew1},
Theorem 1.11, due to S. Bochner and R. S. Phillips (1941), is given, which states that the following construction corresponds to the natural order in the space $ba(X, \Sigma)$ (and, therefore, in the subspaces $ca(X, \Sigma)$, $pfa(X, \Sigma)$):
$$
(\mu_{1} \land \mu_{2})(E) = inf \{ \mu_{1}, \mu_{2} \}(E) = inf ( \mu_{1}(C) + \mu_{2}(E \cap C^{\prime})), 
C^{\prime}=X\setminus C, 
$$ 
where the infimum is taken over all sets $C \subset E, C \in \Sigma$;
$$
(\mu_{1} \lor \mu_{2}) =sup \{ \mu_{1}, \mu_{2} \} = -((-\mu_{1}) \land (-\mu_{2})).  
$$

In this case, $\mu_{1} \land \mu_{2}, \mu_{1} \lor \mu_{2} \in ba(X,\Sigma)$.
%


\begin{definition}
{\it Two positive measures $\mu_{1} , \mu_{2} \in ba(X, \Sigma)$ are called disjoint if $\mu_{1} \wedge \mu_{2}=0$.}
\end{definition}


\begin{definition}
{\it  Two positive measures $\mu_{1} , \mu_{2} \in ba(X, \Sigma)$ are called singular if there are two sets $D_{1}, D_{2}\subset X$,  $\; D_{1}, D_{2}\in \Sigma$}, such that $\mu_{1}(D_{1})=\mu_{1}(X)$, $\; \mu_{2}(D_{2})=\mu_{2}(X)$ and $D_{1}\cap D_{2}=\emptyset$.
\end{definition}

We give in our notation Theorem 1.21 from  ~\cite{YoHew1}.

\begin{theorem}  
 
 1. Let two finitely additive positive measures $\mu_{1},\mu_{2} $ be disjunct. Then for any numbers $\alpha, \beta >0$ there is a set $E\in \Sigma$ such that $\mu_{1}(E) \le \alpha$ and $\mu_{2}(X \setminus E) \le \beta$.
Here you can put $\varepsilon = \alpha = \beta$ and consider $\varepsilon$ "as small as you like".

 2. If at least one of the two finitely additive measures $ \mu_{1}$ and $\mu_{2}$ is not countably additive, then their disjunctness does not imply a singularity.
 
  3. Countably additive measures $\mu_{1}, \mu_{2}$ are disjunct if and only if they are singular.\medskip
\end{theorem}

If the measures $\mu_{1}$ and $\mu_{2}$ are singular, then they are also disjoint (see \cite {YoHew1}).


\begin{definition}
A cycle $K=\{ \mu_{1}, \mu_{2}, ..., \mu_{m} \} $ is called a cycle of disjoint measures if all its cyclic measures are pairwise disjoint,
i.e. $\mu_{i} \wedge \mu_{j} =0$ for all $ i \ne j $.
\end{definition}


\begin{definition}  
Two cycles of measures $K^{1}, K^{2} $ are called disjoint, if each measure from the cycle $ K^{1} $ is disjoint with each measure from the cycle $ K^{2} $.
\end{definition}

If the cycle of disjoint measures $K=\{ \mu_{1}, \mu_{2}, ..., \mu_{m} \} $ is countably additive, then all its cyclic measures are pairwise singular and have pairwise disjoint supports (sets of full measure) $D_{1}, D_{2}, ..., D_{m} \in \Sigma$, that is, $\mu_{i}(D_{i})  $, $ i = 1, ..., m$,  and  $D_{i}\cap D_{j} =\emptyset $ for $i \ne j$.\\

Now we give one of the versions of the classical definition of cycles of a MC in the state space (see, for example, \cite{Rev1}).

\begin{definition}  
A cycle in the set of states $(X, \Sigma)$ of an arbitrary Markov chain is a system of pairwise distinct sets $S =\{ D_{1}, D_{2}, ..., D_{m} \}, D_{i} \in \Sigma$, $ i = 1, ..., m$, such that $p(x, D_{i+1})=1$ for all $x \in D_{i}$, $ i = 1, ..., m -1 $, and $ p (x, D_{1}) = 1 $ for all $ x \in D_{m} $. 
The sets $ D_{i} $, $ i = 1, ..., m $, are called cyclic, and the number $ m $ is the period of the cycle $ S $.

If the sets $ D_{i} $ do not intersect pairwise, that is, $ D_{i} \cap D_{j} = \emptyset $ for $ i \ne j $, then such a cycle is called singular.
\end{definition}

\begin{remark}
Many papers study Markov chains on a measurable set together with some predetermined fixed measure $ \mu $, that is, on $ (X, \Sigma, \mu) $. Accordingly, Markov operators are also constructed on function spaces $ L_{p} (X, \Sigma, \mu) $. With this approach, a refinement naturally appears in the definition of cycles: $ p (x, D_{i + 1}) = 1 $  almost everywhere relative to the measure $ \mu $ for $ x \in D_{i} $, $ i = 1, ..., m-1$, and $ p (x, D_{2}) = 1 $ almost everywhere. Other measures may be used in such definitions. In such cases, they usually write:  [$\mu$] almost everywhere.
\end{remark}

In the present work, we use another apparatus for studying the MC, and no a priori fixed measure $ \mu $ is assumed when specifying the MC. Therefore, there is usually no need to introduce ``almost surely'' or ``almost everywhere'' annotations in the main statements and formulas. But sometimes we can't do without it, such a case will arise in Lemma 4.1 and in Theorem 4.6 at the end of the article.

At the end of the next section, two theorems 4.5 and 4.6 will be proved that show how the cycles of the sets of states of the MC and the cycles of its measures are related.

If we do not require pairwise disjointness (singularity) of the measures of a countably additive cycle, then new, somewhat unexpected objects may appear in the state space of a MC. Let's give a suitable simple example.
 
 {\bf Example 3.1.}
Let the MC be finite, having exactly three states $X = \{ x_{1}, x_{2}, x_{3} \}$ with transition probabilities:
$$
p(x_{1}, x_{1})=1,  p(x_{2}, x_{3})=1,   p(x_{3}, x_{2})=1.  
$$
This means that the MC has in the state space $ X $ one stationary state $\{x_{1}\}$ (we can say that this is a cycle of period $ m = 1 $) and one cycle $\{x_{2}, x_{3}\}$ of the period $ m = 2 $. Within the framework of the operator approach, it is more convenient for us to translate what has been said into the language of measures as follows.

Let $x\in X$ and $E \subset X (E \in \Sigma = 2^{X})$.
Then $p(x_{1}, E)=\delta_{x_{1}}(E)$, $p(x_{2}, E)=\delta_{x_{3}}(E)$, $p(x_{3}, E)=\delta_{x_{2}}(E)$, where $\delta_{x_{i}}(\cdot)$, $i=1, 2, 3$, are the Dirac measures at the points $x_{1}, x_{2},x_{3}$. For an operator $ A $ such a MC we have: $A\delta_{x_{1}}=\delta_{x_{1}}$, $A\delta_{x_{2}}=\delta_{x_{3}}$, $A\delta_{x_{3}}=\delta_{x_{2}}$, i. e.  the family of measures $K=\{\delta_{x_{2}}, \delta_{x_{3}}\}$ is a cycle according to Definition \ref{D:3.2}, and the cyclic measures $\delta_{x_{2}}$ and $\delta_{x_{3}}$ are singular.

Consider one more family of measures $\tilde{K}=\{ \frac{1}{2}\eta_{1}, \frac{1}{2}\eta_{2} \}$, where $\eta_{1}=\delta_{x_{1}}+\delta_{x_{2}}$, $\eta_{2}=\delta_{x_{1}}+\delta_{x_{3}}$.
Then $A\eta_{1}=A(\delta_{x_{1}}+\delta_{x_{2}})$ = $A\delta_{x_{1}}+A\delta_{x_{2}}=\delta_{x_{1}}+\delta_{x_{3}}=\eta_{2}$ and similarly $A\eta_{2}=\eta_{1}$. Since the measures $\eta_{1}$ and $\eta_{2}$ are different, then by Definition \ref{D:3.2}, the family of measures $ \tilde K $ is also a MC cycle different from $ K $. Moreover, the measures $\eta_{1} $ and $ \eta_{2}$ are not disjoint: $\eta_{1} \wedge \eta_{2}=\delta_{x_{1}}\ne 0$. These measures are not singular: their supports $\{x_{1}, x_{2}\}$ и $\{x_{1}, x_{3}\}$ intersect, i. e. $\{x_{1}, x_{2}\} \cap \{x_{1}, x_{3}\}=\{x_{1}\}\ne \emptyset$. \\


\begin{remark} 

Such cycles with intersecting cyclic sets of states, as in Example 3.1, are usually not considered in the classical theory of MC.
\end{remark} 

However, we believe that the study of intersecting cycles of sets is very useful in general theory.
 Research of such cycles is more productive for us in terms of measure cycles.
In this case, instead of intersecting sets of measures, one should consider cycles of measures that are not disjoint.
Our Theorems 4.1, 4.2, 4.3, and 4.7 (proved in Section 4) do not require pairwise disjointness (or singularity) of cyclic measures in measure cycles.


\section {  Main results } 
\label{S:4}


\begin{theorem} 
Any finitely additive cycle of measures for an arbitrary MC is a linearly independent set in the linear space $ba (X, \Sigma)$.
\end{theorem}


{\underline Proof.} We prove by induction.

Consider first two arbitrary different measures  $\mu_{1}, \mu_{2} \in S_{ba}$ (not necessarily cyclic), for which $\|\mu_{1}\|=\| \mu_{2}\|=1$. 
They are obviously linearly independent.
In particular the cycle $K=\{ \mu_{1}, \mu_{2} \}$ consisting of two different measures from $ S_{ba} $, is linearly independent.

Now let the cycle consist of three pairwise different measures: $K=\{ \mu_{1}, \mu_{2}, \mu_{3} \} \subset S_{ba}$.
As we found out above, any two measures of them are linearly independent.

Suppose that one of these three measures is linearly dependent on the other two, let it be the measure $ \mu_{3} $ (the number is not important here). Then there exist numbers $\alpha_{1}, \alpha_{2}, 0 \leq \alpha_{1}, \alpha_{2} \leq 1, \alpha_{1} + \alpha_{2} = 1$, such that the measure $\mu_{3}$ is uniquely representable as a linear combination $\mu_{3} = \alpha_{1} \mu_{1}+ \alpha_{2} \mu_{2}$.   

Let  $\alpha_{1} = 0$.Then $\alpha_{2} = 1$ and $\mu_{3} =  \mu_{2} $ which contradicts the pairwise difference of the three measures. Similarly for $ \alpha_{1} = 1 $. Therefore, we can assume that $ 0 <\alpha_{1}, \alpha_{2} < 1$.

By cycle conditions 
$$
\mu_{1} = A\mu_{3} = A(\alpha_{1} \mu_{1}+ \alpha_{2} \mu_{2}) = \alpha_{1} A\mu_{1}+ \alpha_{2} A\mu_{2} = \alpha_{1} \mu_{2}+ \alpha_{2} \mu_{3}.
$$
Since $\alpha_{2} \ne 0$ from this we get $\mu_{3} = \frac{1}{\alpha_{2}}\mu_{1} - \frac{\alpha_{1}}{\alpha_{2}}\mu_{2}$. Since the decomposition of $ \mu_{3} $ is unique, we have $\alpha_{1} = \frac{1}{\alpha_{2}}, \alpha_{2} = - \frac{\alpha_{1}}{\alpha_{2}} < 0$. Since $\alpha_{1}$ and $ \alpha_{2} $ are positive, we obtain a contradiction in the second equality. Therefore, all three measures $\mu_{1}, \mu_{2}$ and $ \mu_{3} $ are linearly independent.

We turn to the general case.

Let be a cycle of measures $K=\{ \mu_{1}, \mu_{2},..., \mu_{m} \}$ with an arbitrary period $m \geq 3$. We assume that the sets of any $ m-1 $ pieces of measures $ \mu_{i}$ from $K$ are linearly independent. Assume that the measure $\mu_{m}$ (the number is not important) depends linearly on the measures $\mu_{1}, \mu_{2},..., \mu_{m-1}$. Then the measure $ \mu_{m}$ is uniquely represented as
$$ \mu_{m} = \sum_{i=1}^{m-1}\alpha_{i}\mu_{i}, $$
where $0 \le \alpha_{i} \leq 1$ for $i =1, 2,..., m$, $\sum_{i=1}^{m-1}\alpha_{i} = 1$.

Assume that for some $t \in \{ 1, 2, ..., m-1\}$, $\alpha_{t} = 0$ is executed. Then the measure $ \mu_{m}$ is linearly expressed in terms of $m-2 $ pieces of measures $\mu_{i}$, all of them together with $\mu_{m}$ will be $m-1$ piece. This contradicts the assumption that the sets of any $m-1$ pieces of measures $\mu_{i}$ from $K$ are linearly independent. Therefore, all $ \alpha_{i} > 0,\;\; i = 1, 2,..., m-1.$

Now let $ t \in \{ 1, 2, ..., m-1\} $ be $\alpha_{t} = 1$. Then all other $ \alpha_{i} = 0 \; (i \ne t)$ and $\mu_{m} = \alpha_{t} \cdot \mu_{t} = \mu_{t}$, which contradicts the condition of pairwise difference of all measures from the cycle.

So, for all coefficients in the linear decomposition of the measure $ \mu_{m}$ we have $0 < \alpha_{i} < 1,\;\; i = 1, 2,..., m-1$.
 
We apply the operator $ A $ to this decomposition of the measure $ \mu_{m} $ and obtain:

$$
\mu_{1} = A\mu_{m} = \sum_{i=1}^{m-1}\alpha_{i} A\mu_{i} = \sum_{i=1}^{m-1}\alpha_{i} \mu_{i+1} = \alpha_{1} \mu_{2} + \alpha_{2} \mu_{3} + ... + \alpha_{m-1} \mu_{m}.
$$
Therefore, we have ($\alpha_{m-1} \ne 0$):
$$
\mu_{m} = \frac{1}{\alpha_{m-1}} \mu_{1} - \frac{1}{\alpha_{m-1}}\sum_{i=1}^{m-2}\alpha_{i}\mu_{i+1}.
$$

Since the representation for the measure $\mu_{m}$ is unique, here and above we obtain the following relations for the coefficients of the measure $\mu_{2}$:
$$
0 < \alpha_{2} = - \frac{\alpha_{1}}{\alpha_{m-1}} < 0.
$$

It follows from the contradiction obtained that the measure $ \mu_ {m} $ is linearly independent of the other measures of the cycle. Consequently, any other measure $ \mu_ {i} \in K $ is linearly independent of the other measures of the cycle $ K $. The theorem is proved. \\


\begin{theorem} 
Let $K=\{ \mu_{1}, \mu_{2}, ..., \mu_{m} \}$ be a finitely additive cycle of measures for an arbitrary MC. If at least one cyclic measure $\mu_{i}$ is countably additive, then all other cyclic measures in $K$ and their mean measures are also countably additive. Such cycles will be called countably additive.
\end{theorem}



\underline{Proof.}

Since $\mu_{i+1}=A\mu_{i}, i=1, 2, ..., m-1$ and $\mu_{1}=A\mu_{m}$ then the statement of the theorem follows from the fact that the operator $ A $ has the space $ ca (X, \Sigma) $ as its invariant subspace in $ ba (X, \Sigma) $, that is, transforms countably additive measures into countably additive ones.
The countable additivity of the mean measure follows from the fact that $ ca (X, \Sigma) $ is a linear space, i. e. the sum of countably additive measures is also countably additive and a countably additive measure multiplied by a number is also countably additive.

The theorem is proved. \\


\begin{proposition}
There exist  classical Markov chains with purely finitely additive cycles of measures with period $ m \geq 2 $.
\end{proposition}


{\bf Example 4.1.} An example of a classical MC is constructed, for which the existence of a purely finitely additive cycle of measures is proved.


For simplicity, we take a deterministic MC generated by a point transformation.


Let $ X = (0, 1) \cup (1, 2)$, $\Sigma=\mathfrak{B}_{X}$ (Borel $ \sigma $-algebra on $ X $). 
Denote $D_1=(0,1)$, $D_2=(1,2)$.
 Then $D_1\cup D_2 =X $, $D_1\cap D_2 =\emptyset$..

Let's define the transition function of the Markov chain according to the rules:

$p(x,\{1+x^2\})=1$, if $x \in (0, 1)$;  

$p(y,\{(y-1)^2\})=1$, if $y \in (1, 2)$.

Then
$p(x, D_2) = 1$, if $x \in D_1$;
$p(x, D_1) = 1$, if $x \in D_2$.

Therefore, the sets of states $ D_1 $ and $ D_2 $ are cyclic and form a singular cycle           $ S = \{D_{1}, D_{2} \} $ with period $ m = 2 $.

The phase portrait of such a MC is given in Figure 1.

\begin{figure}[h!]
\center{\includegraphics[width=0.7\linewidth]{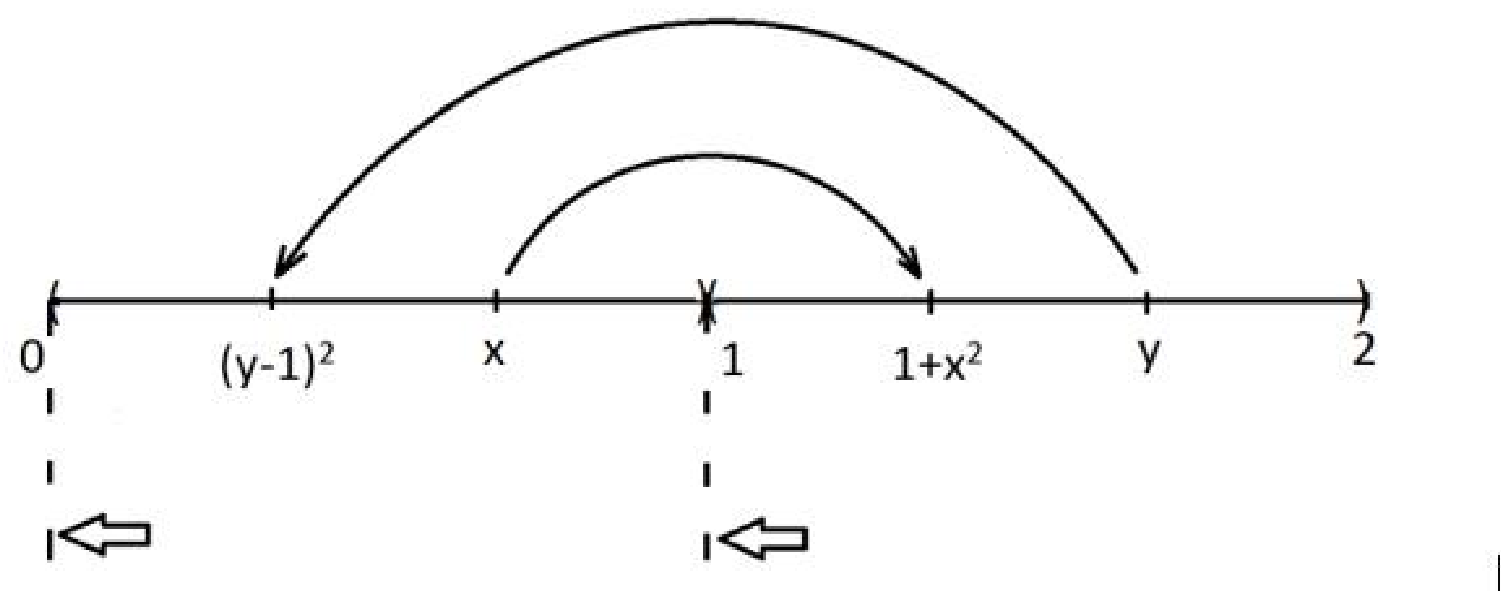}}
\caption{Phase portrait for MC with purely finitely additive cycles.}
\end{figure}

Note that for any trajectory of the Markov chain beginning at the point $ x_0 \in (0,1) $, its subsequence with even numbers tends to zero from the right:
$$
1+x_{0}^{2}, \; 1+x_{0}^{16}, 1+x_{0}^{64}, ... \to 1 ,
$$
and the subsequence with odd numbers tends to one from the right:
$$
x_{0}, \; x_{0}^{4}, x_{0}^{32}, ... \to 0 
$$
(and vice versa, for $x_0 \in (1,2)$).       

By \v{S}idak's Theorem (\cite[Theorem 2.2]{Si1})
for a given MC there exists an invariant finitely additive measure $ \mu = A  \mu \in S_{ba} $. It can be shown that for her $ \mu (D_1) = \mu (D_2) = \frac {1}{2}> 0 $.

We construct two new measures $ \mu_1 $ and $ \mu_2 $ as the restriction of the measure $ \mu $ to the sets $ D_1 $ and $ D_2 $: 
$\mu_1 (E) = \mu(E \cap D_1)$, $\mu_2 (E) = \mu(E \cap D_2)$ for all $ E \subset X $, $ E \in \Sigma $, and $ \mu = \mu_1 + \mu_2 $. The measures $ \mu_1 $ and $ \mu_2 $ are singular and have supports $ D_1 $ and $ D_2 $. It can be proved that $ A \mu_1 = \mu_2 $ and $ A \mu_2 = \mu_1 $. This means that the measures $ \mu_1 $ and $ \mu_2 $ form a disjoint cycle of finitely additive measures $ K = \{\mu_1, \mu_2 \} $.

Let $ 0 <\varepsilon <1 $ and $ D_1^\varepsilon = (0, \varepsilon) $, $ D_2^\varepsilon = (1, 1+ \varepsilon) $. We can get that for any $ \varepsilon $, $ \mu_1 (D_1^\varepsilon) = 1/2 $, $ \mu_2 (D_2^\varepsilon) = 1/2 $. This means that the measures $ \mu_1 $ and $ \mu_2 $ and their mean measure are purely finitely additive. The constructed MC has no invariant countably additive measures.

It can be shown that the singular sets $D_1^\varepsilon$ and $D_2^\varepsilon$ for any $ \varepsilon $ form a cycle of states $S^\varepsilon =\{D_1^\varepsilon, D_2^\varepsilon\}$ and are also supports of measures $ \mu_1 $ and $ \mu_2 $.

It can be proved that the family of all pairwise disjoint invariant finitely additive measures of a given MC has cardinality at least a continuum, i. e. $2^{\aleph_{0}}$. Consequently, the family of all generated pairwise disjoint cycles of purely finitely additive measures of a given MC also has cardinality at least a continuum.

This implies that the algebraic dimension of the set $ \Delta_{ba} $ for a given MC, that is, the cardinality of its Hamel basis, is not less than a continuum.

Note that the family of all purely finitely additive measures `` near '' zero has the cardinality of the hypercontinuum, i. e. is equal to $2^{2^{\aleph_{0}}}$.\\

 
Let us modify the considered MC - add the points $ 0 $ and $ 1 $ to $ X = (0, 1) \cup (1, 2) $ and get $ X = [0, 2) $. Let us determine the possible transitions from these points using the same formulas as the original MC. We get:
 
 $ p (0, \{1 \}) = 1 $,
 $ p (1, \{0 \}) = 1. $
 
 This means that the family of state sets $ S_{0} = \{ \{0\}, \{1\} \} $ for the new MC is a new singular cycle of dimension $ m = 2 $.

It corresponds to a new singular cycle of countably additive measures $ K_{0} = \{\delta_0, \delta_1 \} $, where $ \delta_0 $ and $ \delta_1 $ are Dirac measures at the points $ 0 $ and $ 1 $, respectively. Their mean measure $\mu = \frac{1}{2} \delta_0 + \frac{1}{2} \delta_1$ is countably additive and is the only invariant measure of the modified Markov chain in the class of countably additive measures.

The phase portrait of the new MC is shown in Figure 2.

\begin{figure}[h!]
\center{\includegraphics[width=0.7\linewidth]{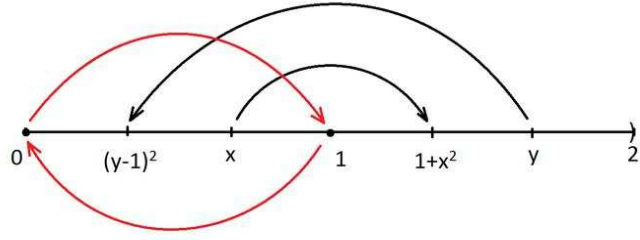}}
\caption{Phase portrait for MC with purely finitely additive cycles
and with one countably additive cycle}
\end{figure}

Note that the whole infinite family of disjoint purely finitely additive cycles of measures $K$ considered above on $X=(0, 1)\cup (1, 2)$ remains the same for the new MC.


\begin{theorem} 
Let $K=\{ \mu_{1}, \mu_{2}, ..., \mu_{m} \} $be a finitely additive cycle of measures for an arbitrary MC. If at least one cyclic measure $\mu_{i}$ is purely finitely additive, then all other cyclic measures in $K$ and their mean measure are also purely finitely additive. Such cycles will be called purely finitely additive.
\end{theorem}

\underline{Proof.}
Let the cyclic measure $\mu_{1}$ be purely finitely additive (the number is not important here) and $\mu_{2} = A\mu_{1}$.

Suppose that the measure $\mu_{2}$ is not purely finitely additive. We decompose the measure $\mu_{2}$ into two components $\mu_{2} = \lambda_{ca} + \lambda_{pfa}$ where $ \lambda_{ca} $ is a countably additive measure, and $ \lambda_{pfa} $ is purely finitely additive. By assumption $\mu_{2} \ne \lambda_{pfa}$ whence $ \lambda_{ca} \ne 0$, $ \lambda_{ca} \geq 0, \lambda_{ca} (X) = \gamma > 0$. 

We apply the operator $ A $ to the measure $ \mu_{2} $
$$
A\mu_{2} = A\lambda_{ca} + A\lambda_{pfa} = \mu_{3}.
$$

The operator $ A $ takes countably additive measures to the same ones and is isometric in the cone of positive measures. It follows from this that the measure $ A \lambda_{ca} $ is countably additive, positive, and $ A \lambda_{ca} (X) = \gamma> 0 $. This means that the measure $\mu_{3}$ also has a positive countably additive component $A\lambda_{ca}$. 


Continuing this procedure further at the last step we get the decomposition
$$
\mu_{1} = A\mu_{m} = A^{m-2}\lambda_{ca} + A^{m-2}\lambda_{pfa},
$$
where the measure $ A^{m-2} \lambda_{ca} $ is countably additive, positive, and $A^{m-2}\lambda_{ca} (X) = \gamma > 0$. 

Thus, the initial measure $\mu_{1}$ has a nonzero countably additive component and, thus, is not purely finitely additive, which contradicts the conditions of the theorem. Therefore, the measure $\mu_{2}$ is also purely finitely additive.

Repeating this procedure sequentially for all the other cyclic measures $\mu_{3}, \mu_{4}, ..., \mu_{m}$ we get that they are all purely finitely additive.

It remains to prove that the mean cyclic measure
$$
\mu = \frac{1}{m}\sum_{i=1}^{m}\mu_{i}
$$ 
is also purely finitely additive. But this follows from the fact that the space of purely finitely additive measures $pfa(X, \Sigma)$ is also linear, which is proved nontrivially in [2] (Theorem 1.17).
The theorem is proved. \\

Now let us present an extended cyclic analogue of the Alexandroff-Yosida-Hewitt decomposition given in Section \ref{S:2}.




\begin{theorem} 
Let $K=\{ \mu_{1}, \mu_{2}, ..., \mu_{m} \}$ be a finitely additive cycle of measures of pairwise disjoint measures with period $ m $ of an arbitrary MC and $\mu_{i} = \mu^{ca}_{i} + \mu^{pfa}_{i}$ a decomposition of cyclic measures into a countably additive component $\mu^{ca}_{i}$ and a purely finitely additive component $\mu^{pfa}_{i}, \; i=1, 2, ..., m$. Then these components are also cyclic, form the cycles $K^{ca}$ and $K^{pfa}$, the cycle $K$ is the coordinate sum of these cycles $K = K^{ca} + K^{pfa}$, and the mean measure of the cycle $K$ is uniquely representable as the sum of its countably additive and purely finitely additive components, which coincide with the mean measures of the cycles $K^{ca}$ and $K^{pfa}$, respectively. Moreover, the cycles $ K ^ {ca} $ and $ K ^ {pfa} $ consist of pairwise disjoint measures and are disjoint with each other, i. e. every measure from $ K^{ca} $ is disjoint with every measure from $ K^{pfa} $.
\end{theorem}

\underline{Proof.}



We denote tuples of countably additive and purely finitely additive components of cyclic measures of a cycle $K$ by the symbols $K^{ca} = \{\mu_1^{ca},\mu_2^{ca}, ..., \mu_m^{ca}\}$ и $K^{pfa} = \{\mu_1^{pfa},\mu_2^{pfa}, ..., \mu_m^{pfa} \}$.
The coordinate-wise sum of these two tuples gives the original cycle $K = K^{ca} + K^{pfa}$. Now we need to show that the measures $\mu_i^{ca}$ and $\mu_i^{pfa}$ are cyclic, that is, the tuples $K^{ca}$ and $K^{pfa}$ are cycles.


Let us prove the theorem step by step.


Assume that some of the measures $\mu^{ca}_{i}$ is zero. Then $\mu_{i}=\mu^{pfa}_{i}$, and according to Theorem 4.3 all other measures
$\mu_{j}=\mu^{pfa}_{j}$, i. e., the cycle  $K=K^{pfa}$, and the theorem is proved. Similarly, for
$\mu^{pfa}_{i}=0$, the cycle $K$ is countably additive by Theorem 4.2, $K=K^{ca}$, and the present theorem is proved.
The main case remains when all
$\mu^{ca}_{i}\ne 0$  and all  $\mu^{pfa}_{i}\ne 0$, which is what we assume below.

Take two arbitrary measures $\mu^{ca}_{i}$ and $\mu^{ca}_{j}$ ($i\ne j$) from $ K ^ {ca} $.
 
 Then
$$0\le \mu^{ca}_{i}  \wedge \mu^{ca}_{j} \le (\mu^{ca}_{i}  + \mu^{pfa}_{i}) \wedge (\mu^{ca}_{j}  + \mu^{pfa}_{j}) =\mu_{i} \wedge \mu_{j}.$$
By the conditions of the theorem, all measures from
$K$ are pairwise disjoint. Therefore, $\mu_{i} \wedge \mu_{j} =0$ and
$\mu^{ca}_{i}  \wedge \mu^{ca}_{j} =0$,  i. e., all measures from 
$K^{ca}$ are pairwise disjoint. And since, as we now assume, all measures from the tuple
$ K^{ca} $ are nonzero, then they are all pairwise distinct.


Similarly, we obtain that all measures from the tuple $ K^{pfa} $ are also pairwise disjoint and distinct.







We emphasize that the tuples of measures $ K^{ca} $ and $ K^{pfa} $ have dimensions $ m $, which coincides with the period $ m $ of the original cycle $ K $.

By the conditions of the theorem, the cycle $ K $ has an (arbitrary) period $ m \in N $. Consequently, each cyclic measure $ \mu_i $ of the cycle $ K $ is an invariant measure of the operator $ A^m $, that is, $ \mu_i = A^{m} \mu_i $, $ i = 1, 2, ..., m $. Take the first cyclic measure with its Alexandrov-Yosida-Hewitt decomposition \cite{YoHew1} $ \mu_1 = \mu_1^ {ca} + \mu_1^ {pfa} $. By \v{S}hidak's Theorem (\cite[Theorem 2.5]{Si1}) both components of the measure $ \mu_1 $ are also invariant measures for the operator $ A^m $, that is, $ \mu_1^{ca} = A^{m} \mu_1^{ca} $, $ \mu_1^{pfa} = A^{m} \mu_1^{pfa} $.

Each of these components generates its own cycle
$$\hat{K}^{ca} = \{\mu_1^{ca}, A \mu_1^{ca}, …, A^{m-1} \mu_1^{ca}\},$$
$$\hat{K}^{pfa} = \{\mu_1^{pfa}, A \mu_1^{pfa}, …, A^{m-1} \mu_1^{pfa} \}.$$

Obviously, the coordinate-wise sum of these two cycles gives the whole cycle $K =\hat{K}^{ca} + \hat{K}^{pfa}$.

Since the measure $\mu_1^{ca}$ is countably additive, then, according to Theorem 4.2, all other cyclic measures of the cycle $ \hat{K}^{ca}$ are countably additive. Since the measure $ \mu_1^{pfa} $ is purely finitely additive, then, according to Theorem 4.3, all other cyclic measures of the cycle $\hat{K}^{pfa}$ are purely finitely additive.









By the uniqueness of the decomposition of any measure into countably additive and purely finitely additive components (see \cite{YoHew1}), we obtain the following equalities (here the symbol $A^{0}$ means the identical operator):

$$\mu_1^{ca} = A^{0} \mu_1^{ca},$$
$$\mu_2^{ca} = A \mu_1^{ca}, ... $$
$$\mu_m^{ca} = A^{m-1} \mu_1^{ca},$$
where on the left are the measures of the tuple $K^{ca}$, and on the right are the cyclic measures of the cycle $ \hat K^{ca}$.

Similar equalities are also true for purely finitely additive components.

From this we get that $K^{ca} =\hat{ K}^{ca}$, $K^{pfa} = \hat{K}^{pfa}$, i.e. tuples $K^{ca}$ и $K^{pfa}$ are cycles, and $K = {K}^{ca} + {K}^{pfa}$. Note that this decomposition of the cycle $ K$ is unique. The main statement of the theorem is proved.

Now the corresponding equalities for the mean measures of cycles are obvious.






In \cite{YoHew1} (Theorem 1.16), it was proved that any countably additive measure is disjoint with any purely finitely additive measure. Therefore, the cycles of the measures $K^{ca}$ and $K^{pfa}$ are disjoint.
Above we showed that all measures from $K^{ca}$ and $K^{pfa}$ are also pairwise disjoint.

The theorem is proved. 


\begin{corollary} 
A finitely additive cycle of measures $ K $ is countably additive if and only if its mean measure is countably additive.
\end{corollary} 


\begin{corollary} 
A finitely additive cycle of measures $ K $ is purely finitely additive if and only if its mean measure is purely finitely additive.
\end{corollary}

Comment.

Under the conditions of Theorem 4.4 just proved, the requirement of pairwise disjointness of cyclic measures in the cycle $ K $ is essential. If we remove it, then the theorem becomes incorrect. 

\begin{theorem} 
Let an arbitrary MC have a singular cycle of the sets of states $S = \{ D_{1}, D_{2}, ..., D_{m} \}$.
Then the MC has a singular cycle of finitely additive measures $K = \{ \mu_1, \mu_2, ..., \mu_m \} \subset S_{ba}$, such that $\mu_i (D_{j})=1$, if $ i=j$, and $ \mu_i (D_{j})=0$, if $ i \ne j $.
\end{theorem}

\underline{Proof.}


We denote the set $D=\cup_{i=1}^{m}D_{i}$. As follows from Definition 3.12, the set $ D $ is stochastically closed for a MC, i. e. $ p (x, D) = 1 $ for all $ x \in D $. Therefore, we can restrict the MC to the set $ D $ by the rule: $ \tilde{p} (x, E) = p (x, E \cap D) $ for all $ x \in D $ and $ E \in \Sigma $. We obtain a new MC on the space $ (D, \Sigma_{D}) $, where the sigma-algebra $ \Sigma_{D} = D \cap \Sigma $.
The transition function $ \tilde{p} (x, E) $ generates the Markov operator $A_{_D}$, which coincides with the original operator $ A $ on the measure space $ba(D, \Sigma_{D})$.





By \v{S}idak ~\cite[theorem 2.2]{Si1}, a restricted MC has in its $ S_{ba} $ invariant finitely additive measure $\mu=A_{_D}\mu$.
Let us expand the invariant measure $ \mu $ into the sum of measures $ \mu_{i} $ constructed in a similar way on cyclic sets $D_{i}$, $i=1, 2, ..., m$, according to the rule: $\mu_{i}(E)=\mu (E\cap D_{i})$ for all $E\in \Sigma_{D}$, and $ i = 1, 2, ..., m $. Obviously, $\mu=\frac{1}{m} \sum_{i=1}^{m}\mu_{i}$ on the $(D, \Sigma_{D})$. 

Since, by the conditions of the theorem, the cyclic sets $ D_{i} $ are pairwise singular, the measures $ \mu_{i} $ are also pairwise singular, the sets $ D_{i} $ are carriers of the measures $ \mu_{i} $ , $ i = 1, 2, ..., m $.

Now we need to show that the measures $ \mu_{i} $, $ i = 1, 2, ..., m $ form a cycle according to Definition \ref{D:3.2}.

Let $E \subset D_{2}, E \in \Sigma_{D}$. The mean measure $ \mu $ is invariant by construction. Therefore,
$$
\mu(E)=A_{_D}\mu(E)=A_{_D}(\frac{1}{m} \sum_{i=1}^{m}\mu_{i})(E)=\frac{1}{m} \sum_{i=1}^{m} A_{_D} \mu_{i}(E)
$$

$$
=\frac{1}{m} \sum_{i=1}^{m} \int\limits_{D} p(x, E) \mu_{i} (dx)
=\frac{1}{m} \sum_{i=1}^{m} \int\limits_{D_{i}} p(x, E) \mu_{i} (dx)
$$
$$
=
\frac{1}{m} \int\limits_{D_{1}} p(x, E) \mu_{1} (dx)=\frac{1}{m} \int\limits_{D} p(x, E) \mu_{1} (dx)=\frac{1}{m} A \mu_{1} (E)
.$$
On the other hand, $$\mu(E)=\frac{1}{m} \sum_{i=1}^{m}\mu_{i}(E)=\frac{1}{m} \mu_{2}(E).$$





Since $\mu = A \mu$, we get $\frac{1}{m} \mu_{2}(E)=\frac{1}{m} A_{_D}\mu_{1}(E)$, from which $\mu_{2}(E)= A_{_D}\mu_{1}(E)$ for all $E \subset D_{2}$, i. e. $\mu_{2} = A_{_D}\mu_{1}$ on the $(D, \Sigma_{D})$.

Similarly, we get $\mu_{i+1}= A_{_D}\mu_{i}$ for $i=2,..., m-1$ and $\mu_{1}= A_{_D}\mu_{m}$ on the $(D, \Sigma_{D})$.

So, the set of measures $K = \{\mu_1, \mu_2,..., \mu_m \}$ forms a cycle of a MC narrowed to the set of States $D=\cup_{i=1}^{m}D_{i}$. 

Now we perform the reverse procedure. Restoring the values of the transition function $p(x, E)$ to the entire measurable space $(X,\Sigma)$.
Continue with $(X,\Sigma)$ constructed measures $\mu$ and $\mu_{i}$, $i=1, 2,..., m$, according to the rule: for all sets $E\subset X\setminus D$, $E\in \Sigma$, we assume $\mu(E)=0$, $\mu(E_{i})=0$, $i=1, 2,..., m$.

Obviously, we get a finite-additive measure cycle defined already on the entire phase space $(X,\Sigma)$.
The theorem is proved.



\begin{remark}  

In Theorem 4.5, we proved that the disjoint cycle $ S $ of the sets of states of a MC generates a disjoint cycle of finitely additive measures $ K $. This statement cannot be strengthened. In the previous Example 4.1, we showed that a MC with a disjoint state cycle $ S $ may not have any invariant (or cyclic) countably additive measure at all, but have invariant (and cyclic) purely finitely additive measures. 
\end{remark}

Now we prove one general Lemma not regardless of Markov chains.

\begin{lemma} 

Let $(X, \Sigma)$ be an arbitrary measurable space on which given:

1) non-negative measurable bounded function $f: X \to R^{1}$ with a condition $0 \leq f(x) \leq 1$ for all $x \in X$;

2) a bounded positive countably additive measure $\mu : \Sigma \to R$ with the condition $0 \leq \mu(E) \leq 1$ for all $E \in \Sigma$ and $\mu(X)=1$ (i.e., a probability measure).

Then, if $$ \int\limits_{X} f(x)\mu(dx)=1, $$ then $f(x)=1 \;\; [\mu] $ almost everywhere, i. e. $\mu \{ x\in X: f(x)=1 \} =1$.
\end{lemma}

\underline{Proof.}

We denote the sets $F=\{ x \in X: f(x)=1\}$, $G=\{ x \in X: f(x)<1\}$. Obviously, $F, G \in \Sigma, F \cap G =\emptyset, X=F\cup G$.
We need to prove that $\mu (F)=1$ or $\mu(G)=0$. Assume the contrary. Let $\mu(G)>0$.

By construction, the set of possible values of the function $ f (x) $ for $ x \in G $ is the half-interval $E=[0, 1)$.
 We split the half-interval $ E $ into a countable sum of pairwise disjoint half-intervals $E_{n}=[1-\frac{1}{n}, 1-\frac{1}{n+1}), n \in N, E=\cup_{n=1}^{\infty}E_{n}$. 
 We take the inverse images of these sets $G_{n}=f^{-1}(E_{n}) = \{ x \in X: f(x) \in E_{n} \} = \{ x \in X: 1-\frac{1}{n} \leq f(x) < 1-\frac{1}{n+1} \} $. The inverse images of $G_{n}, n \in N, $ are also pairwise do not intersect and $G=\cup_{n=1}^{\infty}G_{n}$. Since the measure $\mu$ is countably additive, then $\mu(G)=\mu(\cup_{n=1}^{\infty}G_{n})=\sum_{n=1}^{\infty} \mu (G_{n})$.  By assumption, $\mu(G)>0$. 
 Therefore, there is a number $k \in N$, such that $\mu(G_{k}) = \alpha >0$.
 
 Now let's make integral transformations with ascending inequalities.
$$
1=\int\limits_{X} f(x)\mu(dx)= \int\limits_{X\setminus G_{k}} f(x)\mu(dx) +\int\limits_{G_{k}} f(x)\mu(dx) 
$$
$$
\leq
\int\limits_{X\setminus G_{k}}1 \cdot \mu(dx)+\int\limits_{G_{k}} (1- \frac{1}{k+1})\mu(dx)
$$
$$
=\mu(X\setminus G_{k}) + (1- \frac{1}{k+1})\mu(G_{k})=
\mu(X\setminus G_{k}) + \mu(G_{k}) - \frac{1}{k+1}\mu(G_{k})
$$
$$
=\mu(X) - \frac{1}{k+1} \alpha = 1 - \frac{1}{k+1}\alpha < 1.
$$

Hence we get $ 1 <1 $, i.e. contradiction.

Therefore, $\mu(G)=0$ and $\mu(F)=1$.  The lemma is proved.


If in Lemma 4.1 we abandon the countable additivity of the measure $ \mu $, then the lemma becomes invalid. Let us give a corresponding counterexample.

{\bf Example 4.2.} 
Let a function $f(x)=x$ and a purely finitely additive measure $ \mu $ be given on the interval $X=(0, 1), \Sigma=\mathfrak{B}$, satisfying the condition $\mu((1-\varepsilon, 1))=1$ for any $0< \varepsilon < 1$. 

 Since $1 > f(x)$ for all $x \in (0, 1)$, then for any $0< \varepsilon < 1$  we have
$$
1\geq \int\limits_{(0, 1)} f(x)\mu(dx)= \int\limits_{(0, 1)} x \mu(dx)=\int\limits_{(1-\varepsilon, 1)} x \mu(dx) 
$$
$$
\geq
\int\limits_{(1-\varepsilon, 1)} (1-\varepsilon) \mu(dx) = (1-\varepsilon) \mu((1-\varepsilon), 1) = 1-\varepsilon,
$$
i. e. 
$$1\geq \int\limits_{(0, 1)} f(x)\mu(dx) \geq 1-\varepsilon$$
Hence, $ \int\limits_{(0, 1)} f(x)\mu(dx)=1 $.
However, the set $F=\{ x \in X: f(x)=1\}$ is empty and $\mu (\{ x \in X: f(x)=1\})=0$. 
We can say that the function $f(x)=x$ is not equal to 1 everywhere, not just $ [\mu] $ - almost everywhere on the interval $X=(0, 1)$. 
Lemma 4.1 does not hold for such a pair $f$ and $\mu$.

 
 \begin{theorem}  

Let an arbitrary MC have a disjoint countably additive cycle of measures $K = \{ \mu_1, \mu_2, ..., \mu_m \} \subset S_{ca}$.
Then the MC has a singular cycle of state sets $ S = \{ D_{1}, D_{2}, ..., D_{m} \} $, such that
$\mu_{i}(D_{j})=1$, if $i=j$;  $\mu_{i}(D_{j})=0$, if $i\ne j$,
and $p(x, D_{i+1})=1$ $[\mu_{i}]$ almost everywhere for $x \in D_{i}$, for $i=1, 2, ..., m$, and $p(x, D_1)=1$ $[\mu_{m}]$ almost everywhere for $x \in D_{m}$.

\end{theorem}



\underline{Proof.}
As noted in Theorem 3.1 (with reference to \cite{YoHew1}), if a cycle of countably additive measures is disjoint, then it is also singular. Therefore, there exist pairwise singular (disjoint) supports of cyclic measures $S = \{ D_{1}, D_{2}, ..., D_{m} \}$.  Moreover, $\mu_{i}(D_{j})=1$, if $ i = j $, and $\mu_{i}(D_{j})=0 $, if $ i \ne j $. It is necessary to prove that the family of sets $ S $ also form a cycle by Definition 3.12, and possibly in view of Remark 3.1.


Take the first two adjacent cyclic measures $\mu_{1}, \mu_{2}$ from the cycle $K$. 
Since $\mu_{2}=A\mu_{1}$, then the following equalities are true:
$$1=\mu_{2}( D_{2})=A\mu_{1}(D_{2})=\int\limits_{X} p(x, D_{2}) \mu_{1} (dx)=\int\limits_{D_{1}} p(x, D_{2})\mu_{1} (dx).$$

Note that the measure $\mu_{1}$ is countably additive, $\mu_{1}(D_{1})=1$ and $\mu_{1}\in S_{ca}$.
The transition function of the Markov chain satisfies the inequalities $0\leq p(x, D_{2})\leq 1$ for all $x \in D_{1}$.
The conditions of Lemma 4.1 are fulfilled.
Hence,
 $p(x, D_{2})=1$ $ [\mu_{1}]$  almost everywhere for $x \in D_{1}$.
 
Similarly, one can show that all other pairs of neighboring  measures $\mu_{i}, \mu_{i+1}$ of the cycle $ K $ and pairs of neighboring sets $D_{i}, D_{i+1}$ of the family $ S $ satisfy the conditions 
 $p(x, D_{i+1})=1 \;\; [\mu_{i}]$   almost everywhere for $x \in D_{i}$, for $i=1, 2, ..., m-1$ and $p(x, D_1)=1$  $\;\;[\mu_{m}]$ almost everywhere for $x \in D_{m}$. 
  
 This means that the family of sets $S = \{ D_{1}, D_{2}, ..., D_{m} \}$ is a cycle that satisfies the statement of the theorem.
 The theorem is proved.
 
 \begin{corollary} 
 If the conditions of Theorem 4.6 are satisfied, then
 $p(x, D_{i})=0$  $ [\mu_{i}]$   almost everywhere for $x \in D_{i}$, for $i=1, 2, ..., m$.

 \end{corollary}
 
  \begin{remark} 
  The proof of Theorem 4.6 is based on Lemma 4.1, in which the requirement of countable additivity of the measure $\mu$ is essential (see Example 4.2). Accordingly, this requirement for the measures $\mu_{1}, ...,  \mu_{m}$ is also necessary for Theorem 4.6.
  \end{remark}


\begin{theorem} 
Let an arbitrary MC have one finitely additive cycle of measures $K $ of any period and its mean measure $\mu$ is the only invariant finitely additive measure for the operator $A$. Then the cycle $K$ and its mean measure $\mu$ are countably additive.
\end{theorem}


\underline{Proof.}
Consider a cycle of finitely additive measures $K = \{ \mu_1, \mu_2, ..., \mu_m \}$ and its mean measure $\mu = \frac{1}{m} \sum_{i=1}^{m} \mu_i$. In Section 3 shows that the mean measure $\mu$ of the cycle $K$ is invariant for the operator $A$, i.e. $\mu \in \Delta_{ba}$. By the condition of the theorem, this measure is unique in $\Delta_{ba}$, i. e. $\Delta_{ba}=\{ \mu \}$.


In ~\cite[Theorem 8.3]{Zhd01}, it is proved that if a MC has in $ S_{ba} $ a unique invariant measure $ \mu $, i. e. $\Delta_{ba} = \{ \mu \}$, then this measure is countably additive. Therefore, by Theorem 4.2 and Corollary 4.1, the cycle $ K $ is countably additive.


The theorem is proved.

Therefore, it follows (under the above conditions) that there are no ``single'' purely finitely additive cycles.

\section*{Acknowledgments} 

This work was supported by the Russian Foundation for Basic Research (project No. 20-01-00575-a)






\end{document}